\documentstyle{article}
\begin{document}
\begin{center}
{\Large\bf COMPUTATIONAL ASPECTS OF THE\\ 
GRAVITATIONAL INSTABILITY PROBLEM FOR A\\ [0.1cm]
MULTICOMPONENT COSMOLOGICAL\\ [0.2cm] 
MEDIUM}\\ [1.0cm]
{\bf H.J. Haubold}$^1$ and {\bf A.M. Mathai}$^2$\\ [0.1cm] 
$^1$ Office for Outer Space, United Nations, New York, N.Y., USA\\
[0.3cm]
$^2$ Department of Mathematics and Statistics, McGill University,\\
[0.1cm]
Montreal, P.Q., Canada H3A 2K6
\end{center}
\vspace{1cm}
\noindent
Received...
\clearpage
\noindent
{\bf Summary.} The paper presents results for deriving closed-form
analytic
solutions of the
non-relativistic linear perturbation equations, which govern the
evolution of inhomogeneities in a homogeneous spatially flat 
multicomponent cosmological model. Mathematical methods to derive
computable
forms of
the perturbations are outlined.\par
\vspace{1cm}
\clearpage
\noindent
{\Large\bf 1. Introduction}\par
\medskip
This paper deals with solutions of a system of differential
equations describing the evolution of density perturbations in the
Universe. Before being more specific about the differential
equations, we look at previous work setting the physical
environment for deriving such equations.\par
To formulate a quantitative picture of the evolution of
perturbations outside the horizon it is necessary to solve the
perturbed EINSTEIN equations. The density contrast
$\delta(\vec{x})$ is not a gauge invariant quantity. For 
sub-horizon-sized perturbations, $a\lambda=\lambda_{phys}\ll
H^{-1}$,
this fact is of little consequence, as a NEWTONian analysis is
sufficient; $a$ denotes the scale factor and $H$ the
HUBBLE parameter, respectively. For super-horizon-sized
perturbations, $\lambda_{phys}\gg H^{-1}$, the gauge transformation
for $\delta(\vec{x})$ must be taken into account. In this case a
full general relativistic treatment is required. So far two
approaches are available for finally deriving the differential
equation governing the growth or decay of gravitational
instabilities in an expanding Universe. The metric perturbation
approach was invented by LIFSHITZ (1946) with his derivation of the
synchronous equations. This approach was later used by BARDEEN
(1980) to derive the comoving equations. HAWKING (1966) initiated
an alternative approach to follow the evolution of the density
perturbations in the expanding Universe by employing the general
relativistic fluid flow equations. This method was used by OLSON
(1976) to derive equations which are equivalent to the synchronous
equations, and the same method was further developed by LYTH and
MUKHERJEE (1988) to provide the comoving equations.\par
In order to split up the energy density at a given point in
spacetime into an average value plus a perturbation, one must
specify a spacelike hypersurface on which the averaging is to be
performed. The choice of coordinates in order to define the energy
perturbation throughout spacetime is referred to as a choice of
gauge (BARDEEN 1980). BARDEEN proposed to use the comoving gauge to
deal with the general relativistic equations. Another choice is to
use the synchronous gauge (PEEBLES 1980) which, however, is beset
by the problem of introducing arbitrary "extra gauge mode"
solutions. Actually there is an infinity of synchronous gauges.\par
A choice of LIFSHITZ`s or HAWKING`s approach and synchronous or
comoving gauge having been made sets the framework for deriving
differential equations describing the time evolution of the energy
density perturbation. Here again one has the choice to subscribe to
the perfect fluid approximation or the kinetic description (PEEBLES
1980). In the perfect fluid approximation, which is more physically
transparent, the time evolution of the energy density perturbation
is specified by a differential equation involving the pressure
perturbation as an additional unknown. The equation is of second
order in the comoving gauge, assuming that the pressure
perturbation is known, and it has two solutions revealing a growing
and decaying mode (PEEBLES 1980). Since the equation governing the
amplitude $\delta(t)$ is of second order, there are two solutions.
A given perturbation is expressed as a linear combination of
$\delta_+(t)$ and $\delta_-(t)$. At late times, only the projection
onto the growing mode may be important. Physically, the decaying
mode corresponds to a perturbation with initial overdensity and
velocity arranged so that the initial velocity perturbation
eventually "undoes" the density perturbation. In synchronous gauge
the equation is of third or higher order corresponding to the extra
freedom of choosing the respective gauge.\par
Many efforts have been made to derive approximate solutions of the
equations of relativistic and non-relativistic perturbation theory,
usually valid for early time large-scale and for late time small-
scale inhomogeneities, respectively. the primary purpose of this
paper is to derive closed-form analytic solutions of the non-
relativistic linear perturbation equations, which govern the
evolution of inhomogeneities in homogeneous spatially flat
multicomponent cosmological models. The mathematical method
employed here can also be used to tackle the relativistic version
of the equations involved, which will be done in a subsequent
paper. These closed-form solutions are valid for irregularities on
scales smaller than the horizon. They may be used for analytical
interpolution between known expressions for short and long
wavelength perturbations. In Section 2 we present the physical
parameters relevant to the problem and the fundamental differential
equations for arbitrary polytropic index $\gamma_i$ and a radiation
or matter dominated Universe, respectively. Solutions of this
fundamental equation are provided in all subsequent Sections 3,4,
and 5 for relevant parameter sets including the polytropic index
$\gamma_i$ and the expansion law index $\eta$. To catalogue these
solutions in closed-form, in all Secitons the theory of MEIJER's G-
function will be employed. The following results are a continuation
of the ones by HAUBOLD, MATHAI, and MUECKET (1991) MATHAI (1989),
and NURGALIEV (1986).\par
\bigskip
\noindent
{\Large\bf 2. Some Parameter of Physical Significance}\par
\medskip
The growth of density inhomogeneities can begin as soon as the
Universe is matter-dominated. This would be also the case in a
Universe dominated by non-interacting relic WIMPs. Baryonic
inhomogeneities cannot begin to grow until after decoupling because
until then, baryons are tightly coupled to the photons. After
decoupling, when the Universe is matter-dominated and baryons are
free of the pressure support provided by photons, density
inhomogeneities in the baryons and any other matter components can
grow. Actually the time of matter-radiation equality is the initial
epoch for structure formation. In order to fill in the details of
structure formation one needs "initial data" for that epoch. The
initial data required include, among others, (1) the total amount
of non-relativistic matter in the Universe, quantified by
$\Omega_0$, and (2) the composition of the Universe, as quantified
by the fraction of critical density, $\Omega_i=\rho_i/\rho_c$,
contributed by various components of primordial density
perturbations (i= baryons, WIMPs, relativistic particles, etc.).
Here the critical density $\rho_c$ is the total matter density
of the EINSTEIN-DE SITTER Universe.
Speculations about the earliest history of the Universe have
provided hints as to the appropriate initial data: $\Omega_0=1$
from inflation; $0.015\leq \Omega_B\leq0.15$ and
$\Omega_{WIMP}\sim0.9$ from inflation, primordial nucleosynthesis,
and dynamical arguments. In the following we assume that
$\sum \Omega_i=1, \Omega_i=const.$. The cosmological
medium is
considered to be unbounded and that there may be a uniform
background
of relativistic matter. It is further assumed that the matter is
only slightly perturbed from the background cosmological model.
This assumption may give a good description of the behavior of
matter on large scales even when there may be strongly nonlinear
clustering on small scales (primordial objects). Also assumed is
that matter can be approximated as an ideal fluid with pressure a
function of density alone. It consists of $i$ components having the
densities $\rho_i$ and velocities of sound
$\beta_i^2=dP_i/d\rho_i\propto \rho^{\gamma_i-1}$, when an equation
of state 
\begin{equation}
P_i\propto\rho^{\gamma_i}
\end{equation}
has been taken into account. The $i$ components of the medium are
interrelated through NEWTON`s field equation $\Delta\phi=4\pi
G\sum\rho_i$, containing the combined density $\sum
\rho_i$ of all
components. Superimposed upon an initially homogeneous and
stationary mass distribution shall be a small perturbation,
represented by a sum of plane waves
\begin{equation}
\delta=\frac{\delta\rho_i}{\rho}=\delta_i(t)e^{ikx},
\end{equation}
where $\lambda=2\pi a/k$ defines the wave number $k$ and $\lambda$
is the proper wavelength. In the linear approach the system of
second order differential equations describing the evolution of the
perturbation in the non-relativistic component $i$ is
\begin{equation}
\frac{d^2\delta_i}{dt^2}+2(\frac{\dot{a}}{a})\frac{d\delta_i}{dt}
+k^
2\beta_i^2\delta_i=4\pi G\sum_{j=1}^m\rho_j\delta_j, i=1,...,m,
\end{equation}
where for all perturbations the same wave number k is used. This
equation is valid for all sub-horizon-sized perturbations in any
non-relativistic species, so long as the usual FRIEDMANN equation
for the expansion rate is used. Although the following
investigation of closed-form solutions of the fundamental
differential equation governing the evolution of inhomogeneities in
a multicomponent cosmological medium is quite general, it will be
presented within the context of the inflationary scenario.
Therefore in the following we consider an EINSTEIN-DE SITTER
Universe with zero cosmological constant. After inflation the
FRIEDMANN equation for the cosmological evolution reduces to
\begin{equation}
(\frac{\dot{a}}{a})^2=H^2=\frac{8\pi G}{3}\rho.
\end{equation}
The continuity equation gives $\rho\propto a^{-3}$ and from eq. (4)
one has $a\propto t^{2/3}$ and $\rho_i=\Omega_i/6\pi Gt^2$. During
the expansion the HUBBLE parameter changes as $H=\eta t^{-1}$,
where the expansion law index is
$\eta=\frac{1}{2}$ in the radiation-dominated epoch and
$\eta=\frac{2}{3}$ in the matter-dominated epoch respectively. The
wave number
is
proportional to $a^{-1}$ so that $k^2\beta_i^2=k_i^2t^{2(1-\eta-
\gamma_i)}$ in eq. (3), where now the constants $k_i$ come from
both the wave
vector and the velocity of sound. If $k_i$ = 0 the sound velocity
$\beta_i$ equals zero, in which case the adiabatic index $\gamma_i$
loses its sense. Defining the parameter $\alpha_i=2(2-\eta-
\gamma_i)$ that absorbs the adiabatic index as well as the type of
expansion law we can write for eq. (3) using eq. (4):
\begin{equation}
t^2\ddot{\delta}_i(t)+2\eta
t\dot{\delta}_i(t)+k_i^2t^{\alpha_i}\delta_i=\frac{2}{3}\sum^m_
{j=1}
\Omega_j\delta_j(t),\;\;\;\;\;i=1,\ldots,m,
\end{equation}
and dots denote derivatives with respect to time. We introduce a
time operator $\Delta = t\mbox{d}/\mbox{d}t$ and change the
dependent variable, the density perturbation $\delta_i(t),$ dealing
with the function $\Phi_i$ instead of $\delta_i$ by setting
$\delta_i(t) = t^\alpha\Phi_i(t).$ The equation (5) for the
function $\Phi_i$ is then given by
\[\Delta^2\Phi_i
+ b_i\Phi_i = \frac{2}{3}\sum_{j=1}^m \Omega_j\Phi_j,\] where
\begin{equation}
\Phi_i=t^{-\alpha}\delta_i, \; b_i=k_i^2t^{\alpha_i}-\alpha^2, \;
\alpha_i=2(2-\eta-\gamma_i), \; \alpha=-\left(\frac{2\eta-
1}{2}\right).
\end{equation}
Observing that $\Sigma\Omega_j =1$ and operating on both sides of
the differential equation for the $\Phi_i$ by $\Delta^2$ we obtain
the fundamental equation
\begin{equation}
\Delta^4\Phi_i+\Delta^2(b_i\Phi_i)-
\frac{2}{3}(\Delta^2\Phi_i+b_i\Phi_i)=-\frac{2}{3}\sum_{j=1}^m
b_j\Omega_j\Phi_j, \;\;\;i=1,\ldots m,
\end{equation}
(Mathai, 1989).
As indicated above for $\eta$ the values $\frac{2}{3}$ and
$\frac{1}{2}$
are significant and in what follows it will be assumed that $2\geq
\gamma_i \geq \frac{2}{3}$. We have chosen the range of values of
$\gamma_i$ for physical as well as mathematical reasons as will be
evident later in the analysis of eq. (7).
For some of these parameter values we will consider a
multicomponent medium. Consider the special case
\begin{eqnarray}
b_1 & = & k_1^2t^{\alpha_1}-\frac{(2\eta-1)^2}{4},\\ \nonumber
b_2 & = & b_3=\cdots = b_m=b=k^2t^{\alpha}-\frac{(2\eta-1)^2}{4}
\end{eqnarray}
of the fundamental equation (7). This is the case
considered in Haubold, Mathai and Muecket (1991) in detail. In the
present discussion we will use the same notations as in Haubold,
Mathai and Muecket (1991). 
Let $k=0$ in (8). This is case 4.1 of Haubold, Mathai and Muecket
(1991). For this case the parameters are the following:
\[b_1^*, b_2^*=\pm\left\{ \frac{(2\eta-1)^2}{4\alpha_1^2}\right\}
^{\frac{1}{2}},
b_3^*,b_4^*=\pm\left\{ \frac{1}{\alpha_1^2}\left[\frac{(2\eta-
1)^2}{4}+\frac{2}{3}\right]\right\}^\frac{1}{2},\]
\begin{equation}
a_1^*, a_2^*=-1\pm\left\{ \frac{1}{\alpha_1^2}\left[\frac{(2\eta-
1)^2}{4}+\frac{2}{3}-\frac{2}{3}\Omega_1\right]\right\} ^{\frac{1}
{2}}.
\end{equation}
Table 2.1 gives these parameters for $\eta =\frac{2}{3},
\frac{1}{2}$
and $\gamma _i=\frac{2}{3},1,\frac{4}{3},\frac{5}{3},2.$
\clearpage
\begin{center}
Table 2.1
\end{center}
\medskip
\[\begin{array}{l|l|l|l|l|l} \hline
\eta,\;\gamma_i & \alpha_i & \alpha & b_1^*, b_2^* & b_3^*, b_4^*
&
a_1^*, a_2^*\;for\; \Omega_1=\frac{1}{2}\\ [0.2cm] \hline \hline 
\frac{2}{3}\;\frac{4}{3} & 0 & -\frac{1}{6} & & & \\ [0.2cm] \hline
\frac{2}{3} \; 1 & \frac{2}{3} & -\frac{1}{6} & \pm\frac{1}{4} &
\pm\frac{5}{4} & -1\pm \frac{\sqrt{13}}{4} \\ [0.2cm]\hline 
\frac{2}{3} \; \frac{2}{3} & \frac{4}{3} & -\frac{1}{6} &
\pm\frac{1}{8} & \pm\frac{5}{8} & -1\pm\frac{\sqrt{13}}{8} \\
[0.2cm] \hline 
\frac{2}{3} \; \frac{5}{3} & -\frac{2}{3} & -\frac{1}{6} &
\pm\frac{1}{4} & \pm\frac{5}{4} & -1\pm\frac{\sqrt{13}}{4} \\
[0.2cm] \hline 
\frac{2}{3} \; 2 & -\frac{4}{3} & -\frac{1}{6} & \pm\frac{1}{8}
& \pm\frac{5}{8} & -1\pm\frac{\sqrt{13}}{8} \\ [0.2cm] \hline\hline
\frac{1}{2} \; \frac{4}{3} & \frac{1}{3} & 0 & 0 & \pm\sqrt{6} & -
1\pm\sqrt{3} \\ [0.2cm] \hline 
\frac{1}{2} \; 1 & 1 & 0 & 0 & \pm \sqrt{\frac{2}{3}} & -
1\pm\frac{1}{\sqrt{3}} \\ [0.2cm] \hline 
\frac{1}{2} \; \frac{2}{3} & \frac{5}{3} & 0 & 0 &
\pm\frac{\sqrt{6}}{5} & -1\pm\frac{\sqrt{3}}{5} \\ [0.2cm] \hline 
\frac{1}{2} \; \frac{5}{3} & -\frac{1}{3} & 0 & 0 & \pm\sqrt{6} &
-1\pm\sqrt{3} \\ [0.2cm] \hline 
\frac{1}{2}\;2 & -1 & 0 & 0 & \pm\sqrt{\frac{2}{3}} & -1
\pm\frac{1}{\sqrt{3}}\\ [0.2cm]
\hline \end{array}\]
\medskip
\noindent
When evaluating the solutions explicitly one needs the parameter
differences. These are given in Table 2.2.
\begin{center}
Table 2.2
\end{center}
\[\begin{array}{l|l|l|l|l|l|c|l|l|l|l|l}\hline
\eta & \gamma_i & b_1^* & b_2^* & b_3^* & b_4^* & b_1^*-b_2^* &
b_1^*-
b_3^* & b_1^*-b_4^* & b_2^*-b_3^* & b_2^*-b_4^* & b_3^*-b_4^*\\
[0.2cm] \hline \hline
\frac{2}{3} & 1 & \frac{1}{4} & -\frac{1}{4} & \frac{5}{4} & -
\frac{5}{4} & \frac{1}{2} & -1 & \frac{3}{2} & -\frac{3}{2} & 1 &
\frac{5}{2}\\ [0.2cm]
\frac{2}{3} & \frac{2}{3} & \frac{1}{8} & -\frac{1}{8} &
\frac{5}{8} & -\frac{5}{8} & \frac{1}{4} & -\frac{1}{2} &
\frac{3}{4} & -\frac{3}{4} & \frac{1}{2} & \frac{5}{4}\\ [0.2cm]
\frac{2}{3} & \frac{5}{3} & \frac{1}{4} & -\frac{1}{4} &
\frac{5}{4} & -\frac{5}{4} & \frac{1}{2} & -1 & \frac{3}{2} & -
\frac{3}{2} & 1 & \frac{5}{2}\\ [0.2cm]
\frac{2}{3} & 2 & \frac{1}{8} & -\frac{1}{8} & \frac{5}{8} &
-\frac{5}{8} & \frac{1}{4} & -\frac{1}{2} & \frac{3}{4} &
-\frac{3}{4} &
\frac{1}{2} & \frac{5}{4}\\[0.2cm] \hline \hline  
\frac{1}{2} & \frac{4}{3} & 0 & 0 & \sqrt{6} & -\sqrt{6} & 0 & -
\sqrt{6} & \sqrt{6} & -\sqrt{6} & \sqrt{6} & 2\sqrt{6}\\ [0.2cm]
\frac{1}{2} & 1 & 0 & 0 & \sqrt{\frac{2}{3}} & -\sqrt{\frac{2}{3}}
& 0 & -\sqrt{\frac{2}{3}} & \sqrt\frac{2}{3} & -\sqrt{\frac{2}{3}}
&
\sqrt{\frac{2}{3}} & 2\sqrt{\frac{2}{3}}\\ [0.2cm]
\frac{1}{2} & \frac{2}{3} & 0 & 0 & \frac{\sqrt{6}}{5} & -
\frac{\sqrt{6}}{5} & 0 & -\frac{\sqrt{6}}{5} & \frac{\sqrt{6}}{5}
& -\frac{\sqrt{6}}{5} & \frac{\sqrt{6}}{5} & 2\frac{\sqrt{6}}{5}\\
[0.2cm]
\frac{1}{2} & \frac{5}{3} & 0 & 0 & \sqrt{6} & -\sqrt{6} & 0 & -
\sqrt{6} & \sqrt{6} &-\sqrt{6} & \sqrt{6} & 2\sqrt{6}\\ [0.2cm]
\frac{1}{2} & 2 & 0 & 0 & \sqrt{\frac{2}{3}} & -\sqrt{\frac{2}{3}}
& 0 & -\sqrt{\frac{2}{3}} & \sqrt{\frac{2}{3}} & -
\sqrt{\frac{2}{3}} & \sqrt{\frac{2}{3}} & 2\sqrt{\frac{2}{3}}\\
[0.2cm]\hline
\end{array}\]
{\Large\bf 3. Solution when $\alpha_1=0$}\par
\medskip
Note that when $\eta=\frac{2}{3}$ and $\gamma_i=\frac{4}{3}$ one
can
have $\alpha_i=0$ in equation (6). The exact solution given in
equation (5.3) of
Haubold, Mathai and Muecket (1991) is for the situation
$\alpha_i\neq 0$. Hence this case needs separate discussion. When
$\alpha_1=0$ the fundamental equation (7) reduces to the form
\begin{eqnarray}
\Delta^4\Phi_1 & - &\left[\frac{(2\eta-1)^2}{2}+\frac{2}{3}-
k_1^2\right]\Delta^2\Phi_1+\left[\frac{(2\eta-
1)^4}{16}+\frac{2}{3}\frac{(2\eta-
1)^2}{4}\right.\\ \nonumber
& & + \left.\left(\frac{2}{3}\Omega_1-\frac{(2\eta-1)^2}{4}-
\frac{2}{3}\right)k_1^2\right]\Phi_1=0.
\end{eqnarray}
In this case the general solution is of the form\\
\begin{equation}
\Phi_1=c_1t^{d_1}+c_2t^{d_2}+c_3t^{d_3} +c_4^{d_4},
\end{equation}
where $c_1, c_2, c_3, c_4$ are arbitrary constants and
$d_1,d_2,d_3,d_4$
are the roots of the equation
\begin{eqnarray}
x^4 & - & \left[\frac{(2\eta-1)^2}{2}+\frac{2}{3}-
k_1^2\right]x^2\\ \nonumber
& & +\left[\frac{(2\eta-
1)^4}{16}+\frac{2}{3}\frac{(2\eta-1)^2}{4}+\left\{
\frac{2}{3}\Omega_1-\frac{(2\eta-1)^2}{4}-\frac{2}{3}\right\}
k_1^2\right]=0.
\end{eqnarray}
They can be seen to be the following:
\begin{equation}
d_1,d_2=\pm\left[\frac{(2\eta-1)^2}{4}+\frac{1}{3}-
\frac{k_1^2}{2}+\left\{ (\frac{1}{3}+\frac{k_1^2}{2})-
\frac{2}{3}\Omega_1\right\} ^{\frac{1}{2}}\right]^{\frac{1}{2}};
\end{equation}
\begin{equation}
d_3,d_4=\pm\left[\frac{(2\eta-1)^2}{4}+\frac{1}{3}-
\frac{k_1^2}{2}-\left\{ (\frac{1}{3}+\frac{k_1^2}{2})-
\frac{2}{3}\Omega_1\right\} ^{\frac{1}{2}}\right]^{\frac{1}{2}}.
\end{equation}\par
\bigskip
\noindent
{\Large\bf 4. Solution when the parameters differ by\par
\hspace{0.2cm}integers}\par
\medskip
The general solution given in Haubold, Mathai and Muecket (1991),
equations (5.2) and (5.3), are for finite values of t and for the
cases that the $b_j^*$'s do not differ by integers. From Table 2.2
note that at two points, namely $(\eta=\frac{2}{3},
\gamma_i=\frac{5}{3})$ and $(\eta=\frac{2}{3}, \gamma_i=1)$ one has
$b_1^*-b_3^*=-1$ and $b_4^*-b_2^*=-1.$ This means that $G_1$
and $G_4$ of (5.3) of Haubold, Mathai and Muecket (1991) need
modifications. At all other points the $G_j$'s remain the same.\par
\clearpage
\noindent
{\large\bf 4.1. Modifications of the $G_j's$ for the cases
$(\eta=\frac{2}{3},\gamma_i=1)$}\par 
\hspace{0.3cm}{\large\bf and $(\eta=\frac{2}{3},
\gamma_i=\frac{5}{3})$}\par
\medskip
\noindent
Consider, for $x=k_1^2 t^{\alpha_1}/ \alpha_1^2,
\alpha_1\neq 0,$ 
\begin{eqnarray}
G_1 & = & G_{2,4}^{1,2}\left(x\mid^{a_1^*+1, a_2^*+1}_{b_1^*,
b_2^*, b_3^*, b_4^*}\right)\\ \nonumber
& = & \frac{1}{2\pi i}\int_L\frac{\Gamma(\frac{1}{4}+s)\Gamma(-
a_1^*-s)\Gamma(-a_2^*-s)x^{-
s}ds}{\Gamma(\frac{5}{4}-s)\Gamma(-\frac{1}{4}-
s)\Gamma(\frac{9}{4}-s)}.
\end{eqnarray}
Note that a zero coming from $\Gamma(-\frac{1}{4}-s)$ at $s=-
\frac{1}{4}$ coincides with a pole coming from
$\Gamma(\frac{1}{4}+s)$ at $s=-\frac{1}{4}$. This can be removed by
rewriting as follows:
\begin{equation}
G_1=-\frac{1}{2\pi i}\int_L\frac{\Gamma(\frac{5}{4}+s)\Gamma(-
a_1^*-s)\Gamma(-a_2^*-s)x^{-s}ds}{\Gamma(\frac{5}{4}-
s)\Gamma(\frac{3}{4}-s)\Gamma(\frac{9}{4}-s)}.
\end{equation}
Evaluating (16) as the sum of the residues at the poles of
$\Gamma(\frac{5}{4}+s)$ one has
\begin{eqnarray}
G_1 & = &
-x^{\frac{5}{4}}\frac{\Gamma(-a_1^*+\frac{5}{4})\Gamma(-
a_2^*+\frac{5}{4})}{\Gamma(\frac{10}{4})\Gamma(2)\Gamma(\frac{14}{
4})}\\ \nonumber
& & \times\; _2F_3(-a_1^*+\frac{5}{4}, -a_2^*+\frac{5}{4};
\frac{10}{4}, 2, \frac{14}{4}; -x).
\end{eqnarray}
Evaluating $G_4$ also the same way one has the following:
\begin{eqnarray}
G_4 & = & G_{2,4}^{1,2}\left(x\mid^{a_1^*+1, a_2^*+1}_{b_4^*,
b_1^*, b_2^*, b_3^*}\right)\\ \nonumber
& = & \frac{1}{2\pi i}\int_L\frac{\Gamma(-\frac{5}{4}+s)\Gamma(-
a_1^*-s)\Gamma(-a_2^*-s)}{\Gamma(\frac{3}{4}-s)\Gamma(\frac{5}{4}-
s)\Gamma(-\frac{1}{4}-s)}x^{-s}ds\\ \nonumber
& = & -\frac{1}{2\pi i}\int_L\frac{\Gamma(-\frac{1}{4}+s)\Gamma(-
a_1^*-s)\Gamma(-a_2^*-s)}{\Gamma(\frac{3}{4}-s)\Gamma(\frac{9}{4}-
s)\Gamma(-\frac{1}{4}-s)}x^{-s}ds\\ \nonumber
& = & -x^{-
\frac{1}{4}}\frac{\Gamma(-\frac{1}{4}-a_1^*)\Gamma(-
\frac{1}{4}-a_2^*)}{\Gamma(\frac{1}{2})\Gamma(2)\Gamma(-
\frac{1}{2})}\\ \nonumber
& & \times\; _2F_3(-\frac{1}{4}-a_1^*, -\frac{1}{4}-a_2^*;
\frac{1}{2},2,-\frac{1}{2};-x).
\end{eqnarray}
Thus the complete solution in this case is of the form
\begin{equation}
\Phi_1=c_1G_1+c_2G_2+c_3G_3+c_4G_4,
\end{equation}
where $c_1, c_2, c_3, c_4$ are arbitrary constants, $G_1$ and $G_4$
are given in (17) and (18) respectively and $G_2$ and $G_3$ are
given in equation (5.3) of Haubold, Mathai and Muecket (1991).\par
\clearpage
\noindent
{\Large\bf 5. Solution near $\infty$ for the parameter values
of\par
\hspace{0.2cm}table
2.1}\par
\medskip
Except for the points $(\eta, \gamma_i)=(\frac{2}{3},\frac{2}{3}),
(\frac{2}{3}, 2)$ in all other cases some of the parameters differ
by integers. In five cases of Table 2.2 the poles in the integrand
can be up to order 2 and in two cases the poles can be of order
up to 3. The general solution near $\infty$ is of the form
\begin{equation}
\Phi_1=f_1F_1+f_2F_2+f_3F_3+f_4F_4,
\end{equation}
where $f_1, f_2, f_3, f_4$
are arbitrary constants and $F_j$ are given in the equation (5.6)
of Haubold, Mathai and Muecket (1991). For example,
\begin{equation}
F_1=G_{2,4}^{4,1}\left(x\mid^{1+a_1^*, 1+a_2^*}_{b_1^*,
b_2^*, b_3^*, b_4^*}\right), x=k_1^2t^{\alpha_1}/\alpha_1^2,
\alpha_1\neq0.
\end{equation}
$F_2$ is available from $F_1$ by interchanging $a_1^*$ and $a_2^*$.
$F_3$ and $F_4$ are available from $F_1$ by replacing $x$ by $x
e^{i\pi}$ and $xe^{-i\pi}$ respectively. We will evaluate
$F_1$ for the parameter values which differ by integers.\par
\medskip
\noindent
{\large\bf 5.1. $F_1$ for $(\eta=\frac{1}{2},
\gamma_i=\frac{4}{3})$}\par 
\smallskip
\noindent
In this case one has
\begin{equation}
F_1=\frac{1}{2\pi
i}\int_L\frac{\Gamma(s)\Gamma(s)\Gamma(\sqrt{6}+s)\Gamma(-
\sqrt{6}+s)\Gamma(-a_1^*-s)}{\Gamma(1+a_2^*+s)}x^{-s}ds.
\end{equation}
Note that at $s=0,-1,-2,...$ the integrand has poles of order 2
each. 
All other poles are of order 1 each. Thus
\begin{equation}
F_1=H_1+H_2+H_3,
\end{equation}
where $H_2$ and $H_3$ are the sums of the residues at the poles of
$\Gamma(\sqrt{6}+s)$ and $\Gamma(-\sqrt{6}+s)$ respectively and
$H_1$ is the sum of the residues at the poles of $\Gamma^2(s).$ For
all the cases of simple poles including $H_2$ and $H_3$ one has the
result as follows:
\begin{eqnarray}
R_j & = & x^{b_j^*}\left[\Pi_{k=1}^{'4} \Gamma(b_k^*-
b_j^*)\right]\frac{\Gamma(-a_1^*+b_j^*)}{\Gamma(1+a_2^*-b_j^*)}
\\ \nonumber
& & \times\;_2F_3(-a_1^*+b_j^*,-a_2^*+b_j^*; 1-
b_1^*+b_j^*,\ldots,\#,\ldots,1-b_4^*+b_j^*;-x),
\end{eqnarray}
where $\Pi'$ denotes the absence of the gamma $\Gamma(b_k^*-b_k^*)$
and \# denotes the absence of the element $1-b_j^*-b_j^*$. Thus one
has,
\begin{eqnarray}
H_2 & = & x^{\sqrt{6}}\Gamma^2(-\sqrt{6})\Gamma(-2
\sqrt{6})\frac{\Gamma(-a_1^*+\sqrt{6})}{\Gamma(1+a_2^*-\sqrt{6})}
\\ \nonumber
& & \times
\;_2F_3(-a_1^*+\sqrt{6},-a_2^*+\sqrt{6};1+\sqrt{6},1+\sqrt{6},1+2
\sqrt{6};-x)
\end{eqnarray}
and
\begin{eqnarray}
H_3 & = & x^{-
\sqrt{6}}\Gamma^2(\sqrt{6})\Gamma(2\sqrt{6})\frac{\Gamma(-a_1^*-
\sqrt{6})}{\Gamma(1+a_2^*+\sqrt{6})}\\ \nonumber
& & \times\; _2F_3(-a_1^*-\sqrt{6},-a_2^*-\sqrt{6};1-\sqrt{6},
1-\sqrt{6}, 1-
2\sqrt{6};-x).
\end{eqnarray}
Now consider the evaluation of $H_1$. Note that for all parameter
combinations $(\eta,\gamma_i)=(\frac{1}{2},\frac{4}{3}),
(\frac{1}{2},1),
(\frac{1}{2},\frac{2}{3}),(\frac{1}{2},\frac{5}{3}),
(\frac{1}{2},2)$
we have $b_1^*=0=b_2^*$. Hence we can write $H_1$ in general terms
as follows,
\begin{equation}
H_1=\frac{1}{2\pi
i}\int_L\frac{\Gamma^2(s)\Gamma(b_3^*+s)\Gamma(b_4^*+s)\Gamma(-a_
1^*
-s)}{\Gamma(1+a_2^*+s)}x^{-s}ds.
\end{equation}
The general technique of evaluating such integrals is available in
Mathai and Saxena (1973). The solution is the following
\begin{eqnarray}
H_1 & = & \sum_{\nu=0}^\infty\left[\frac{1}{(\nu!)^2}-
(\mbox{ln}x)\left\{ 2\psi(\nu+1)+\psi(b_3^* -
\nu)+\psi(b_4^*-\nu)\right. \right.\\ \nonumber
& &
\left. \left. +
\psi(-a_1^*+\nu)-\psi(1+a_2^*-\nu)\right\}\right]
\\ \nonumber
& & \times \Gamma(b_3^*-\nu)\Gamma(b_4^*-\nu)\frac{\Gamma(-
a_1^*+\nu)}{\Gamma(1+a_2^*-\nu)}\frac{x^\nu}{(\nu!)^2},
\end{eqnarray}
where the $b_3^*$ and $b_4^*$ are given in Table 2.1 for the
specific parameter combinations and $\psi()$ is a psi function.
Note that the first part in (27) can be written as a $_2F_3(-
a_1^*,-a_2^*;1,1-b_3^*,1-b_4^*;-x).$\par
\medskip
\noindent
{\large\bf 5.2. $F_1$ for $(\eta=\frac{2}{3},
\gamma_i=1)$ and
$(\eta=\frac{2}{3},
\gamma_i=\frac{5}{3})$}\par
\smallskip
\noindent
In these two cases $F_1$ is of the following form:
\begin{eqnarray}
F_1 & = & \frac{1}{2\pi i}\int_L\frac{\Gamma(\frac{1}{4}+s)\Gamma(-
\frac{1}{4}+s)\Gamma(\frac{5}{4}+s)\Gamma(-\frac{5}{4}+s)\Gamma(-
a_1^*-s)x^{-s}ds}{\Gamma(1+a_2^*+s)}\\ \nonumber
& = & \frac{1}{2\pi i}\int_L\frac{\Gamma^2(\frac{5}{4}+s)\Gamma^2(-
\frac{1}{4}+s)\Gamma(-a_1^*-s)}{(s+\frac{1}{4})(s-\frac{5}{4})
\Gamma(1+a_2^*+s)}x^{-s}ds.
\end{eqnarray}
There are poles of order one each at $s=-\frac{1}{4}$ and
$s=\frac{5}{4}$ respectively and poles of order 2 each at $s=-
\frac{5}{4}-\nu, s=\frac{1}{4}-\lambda,
\nu=0,1,\ldots,\lambda=0,1,\ldots$ respectively. Thus $F_1$ is
evalueated as the sum of the residues at all these poles. Let us
denote these residues by $R_1,R_2,R_3$ and $R_4$ respectively. Then
\begin{equation}
F_1=R_1+R_2+R_3+R_4,
\end{equation}
where
\begin{eqnarray}
R_1 & = &\mbox{residue at}\; s=-\frac{1}{4}\\ \nonumber
& = & -\frac{\Gamma^2(-\frac{1}{2})\Gamma(-
a_1^*+\frac{1}{4})}{(\frac{3}{2})\Gamma(\frac{3}{4}+a_2^*)}x^
{\frac{1}{4}}=-\frac{2}{3}\frac{\Gamma^2(-\frac{1}{2})\Gamma(-
a_1^*+\frac{1}{4})}{\Gamma(\frac{3}{4}+a_2^*)}x^{\frac{1}{4}};\\
\nonumber
R_2 & = &\mbox{residue at}\; s=\frac{5}{4}\\ \nonumber
& = & \frac{2}{3}\frac{\Gamma^2(\frac{5}{2})\Gamma(-
a_1^*+\frac{5}{4})}{\Gamma(\frac{9}{4}+a_2^*)}x^{-
\frac{5}{4}};\\ \nonumber
R_3 & = & x^{\frac{5}{4}}\sum_{\nu=0}^\infty\left\{ -
\mbox{ln}x+2\psi (\nu+1)+2\psi(-\frac{3}{2}-\nu)+\psi (-
a_1^*+\frac{5}{4}+\nu )\right.\\ \nonumber
& & \left.-\psi(1+a_2^*-\frac{5}{4}-
\nu)+\frac{1}{\nu+1}+\frac{1}{(\frac{5}{4}+\nu)}\right\}\\
\nonumber
& & \times\frac{\Gamma^2(-\frac{3}{2}-\nu)\Gamma(-
a_1^*+\frac{5}{4}+\nu)}{\Gamma(1+a_2^*-\frac{5}{4}-
\nu)(1+\nu)(\frac{5}{2}+\nu)}\frac{x^\nu}{(\nu!)^2};\\
\nonumber
R_4 & = & x^{-\frac{1}{4}}\sum_{\nu=0}^\infty\left\{ -
\mbox{ln}x+2\psi(\nu+1)+2\psi(\frac{3}{2}-\nu)\right. \\ \nonumber
& & \left.+\psi(-a_1^*-\frac{1}{4}+\nu)-\psi(\frac{5}{4}+a_2^*-
\nu)+\frac{1}{(-\frac{1}{2}+\nu)}+\frac{1}{1+\nu}\right\} \\
\nonumber
& & \times \frac{\Gamma^2(\frac{3}{2}-\nu)\Gamma(-a_1^*-
\frac{1}{4}+\nu)}{(\frac{1}{2}-\nu)(-1-\nu)\Gamma(\frac{5}{4}+a_2
^*-
\nu)}\frac{x^\nu}{(\nu!)^2}.
\end{eqnarray}
As before $F_2$ is available from the $F_1$ of (30) by
interchanging $a_1^*$ and $a_2^*$. $F_3$ and $F_4$ are available
from $F_1$ replacing $x$ by $e^{i\pi}x$ and $e^{-
i\pi}x$ respectively. This completes the evaluation of
$\Phi_1$ for all finite values as well as for values near $\infty$
for all the parameter combinations given in Tables 2.1 and 2.2.\par
\bigskip
\noindent
{\Large\bf 6. Conclusion}\par
\medskip
\noindent
We have presented closed-form solutions of the non-relativistic
linear perturbation equations which govern the evolution of
inhomogeneities in a spatially flat multicomponent cosmological
medium. The general solutions are catalogued according to the
polytropic index $\gamma_i$ and the expansion law index $\eta$ of
the multicomponent medium . All general solutions are expressed in
terms of MEIJER's G-function and their MELLIN-BARNES integral
representation. The proper use of this function simplifies the
derivation of solutions of the linear perturbation equations and
opens ways for its numerical computation. In this regard the paper
leaves room for further work.\par
\bigskip
\begin{center}
Acknowledgments
\end{center}
The authors would like to thank the Natural Sciences and the
Engineering Research Councel of Canada for financial assistance for
this research project.
\clearpage
\noindent
{\Large\bf References}\par
\medskip
\noindent
Bardeen, J.M.: 1980, Phys. Rev., {\bf D22}, 1882\par
\medskip
\noindent
Haubold, H.J., Mathai, A.M., Muecket, J.P.: 1991, Astron.
Nachr. {\bf 312},1\par
\medskip
\noindent
Hawking, S.W.: 1966, Ap.J., {\bf 145}, 544\par
\medskip
\noindent
Lifshitz, E.M.: 1946, J. Phys. (Moscow), {\bf 10}, 116\par
\medskip
\noindent
Lyth, D.H., Mukherjee, M.: 1988, Phys. Rev., {\bf D38}, 485\par
\medskip
\noindent 
Mathai, A.M.: 1989, Studies Appl. Math. {\bf 80},75.\par
\medskip
\noindent
Mathai, A.M., Saxena, R.K.: 1973. Generalized
Hypergeometric Functions\par 
with Applications in Statistics and Physical
Sciences. Lecture Notes in\par 
Mathematics, Vol. 348, Springer-Verlag, Berlin-Heidelberg-New
York.\par
\medskip
\noindent
Nurgaliev, I.S.: 1986, Sov. Astron. Lett. {\bf 12}, 73\par
\medskip
\noindent
Olson, D.W.: 1976, Phys. Rev., {\bf D14}, 327\par
\medskip
\noindent
Peebles, P.J.E.: 1980, The Large-Scale Structure of the
Universe,\par 
Princeton University Press, Princeton\par 
\end{document}